\documentclass[11pt]{elsart}

\usepackage{url}

\usepackage{amssymb}
\usepackage{times}
\usepackage{enumerate}

{\bfseries}{\itshape}
{\bfseries}{\itshape}

\newcommand {\bbox}{\rule{0.6em}{0.6em}}

\newcommand{\boxi}{\ensuremath{{box}}}

\newcommand{\tw}{\ensuremath{{tw}}}
\newcommand{\bag}{\ensuremath{\mathrm{b}}}

\newcommand{\level}{\ensuremath{h}}

\newcommand{\core}{\ensuremath{c}}
\newcommand{\finter}{\ensuremath{f}}
\newcommand{\nr}[4]{\ensuremath{(\{{#1}_i : i \in #2\},\, #3,\, #4)}}
\newcommand{\floor}[1]{\ensuremath{\lfloor #1 \rfloor}}
\newcommand{\ignore}[1]{}

\renewenvironment{pf}
{\begin{oldproof}}{\hfill\bbox\end{oldproof}}

\begin{document}

\begin{frontmatter}

\title{Boxicity and Treewidth}

\author[sun]{L. Sunil Chandran} 
and 
\author[ns]{Naveen Sivadasan} 
\address[sun]{Indian Institute of Science,
Dept. of Computer Science and Automation,
Bangalore  560012, India.  email: \emph{sunil@csa.iisc.ernet.in}}
\address[ns]{Strand Genomics,
237, Sir. C. V. Raman Avenue, Rajmahal Vilas,
Bangalore  560080, India.  email: \emph{naveen@strandgenomics.com}}

\date{}
\bibliographystyle{plain}

\pagestyle{plain}
\pagenumbering{arabic}
\maketitle

\begin {abstract}

An axis-parallel  $b$--dimensional box is a Cartesian product 
$R_1 \times R_2 \times \cdots \times R_b$ where $R_i$ (for $1 \le i \le b$)
is a closed interval of the form $[a_i, b_i]$ on the real line.
For a  graph $G$, its \emph{boxicity} $\boxi(G)$ is the minimum dimension $b$, such that
$G$ is representable as the intersection graph of 
(axis--parallel) boxes in $b$--dimensional space. The concept of
boxicity finds applications in various areas such as ecology, operation research etc.
Though many authors have investigated this concept, not much is known about
the boxicity of many well-known graph classes (except for a couple of
cases) perhaps due to lack of effective approaches. Also,
little is known about the structure imposed on a graph by its high boxicity.

The
concepts of \emph{tree decomposition} and \emph{treewidth} play a very important
role in modern graph theory and has many applications to computer science.
In this paper, we relate the seemingly unrelated concepts of
treewidth and boxicity. Our main result is that, for any graph $G$,
$\boxi(G) \le \tw(G) + 2$, where $\boxi(G)$ and $\tw(G)$ denote the
boxicity and treewidth of $G$ respectively. We also show that this upper
bound is (almost) tight. Since treewidth and tree decompositions
are extensively studied concepts, our result leads to various interesting
consequences, like bounding the boxicity of many well known graph classes,
such as chordal graphs, circular arc graphs, AT-free graphs, co--comparability graphs etc.
 For all these
graph classes, no bounds on their boxicity were known previously.
All our bounds are shown to be tight up to small constant factors.
An algorithmic consequence of our result is a linear time algorithm to construct
a box representation for graphs of bounded treewidth in a space of constant dimension.

Another consequence of our main result is that,
if the boxicity of a graph is $b \ge 3$, then there exists a simple cycle of length
at least $b-3$ as well as an induced cycle of length at least 
$\lfloor\log_{\Delta}(b-2) \rfloor + 2$, where $\Delta$ is its maximum degree.
We also relate boxicity with the cardinality of minimum vertex cover, minimum feedback
vertex cover etc. Another structural consequence is that,
for any fixed planar graph $H$,
there is a constant $c(H)$ such that, if $\boxi(G) \ge c(H)$ then
$H$ is a minor of $G$.

\end {abstract}

\end{frontmatter}

\section{Introduction}

\subsection{Boxicity}
Let ${\mathcal{F}} = \{ S_x \subseteq U : x \in V \}$ be a family of 
subsets of a universe $U$, where $V$ is an index set. The intersection graph 
$\Omega({\mathcal{F}})$ of ${\mathcal{F}}$ 
has $V$ as node set, and two distinct nodes
$x$ and $y$ are adjacent if and only if $S_x \cap S_y \ne \emptyset$.

Representations of graphs as the intersection graphs of various
geometrical objects is a well studied topic in graph theory.  
A prime example of a graph class defined in this way is the class of 
interval graphs: {\it A graph $G$ is an interval graph \emph{if and only if} 
$G$ has an interval realization: i.e., each
node of  $G$ can be associated to an interval on the real 
line such that two intervals intersect if and only if
the corresponding nodes are adjacent.}
Motivated by theoretical as well as practical considerations, 
graph theorists have tried to generalize the concept of interval graphs
in various ways. 
One such generalization is the concept of \emph{boxicity} defined as follows.

An axis-parallel  $b$--dimensional box is a Cartesian product 
$R_1 \times R_2 \times \cdots \times R_b$ where $R_i$ (for $1 \le i \le b$)
is a closed interval of the form $[a_i, b_i]$ on the real line.
For a  graph $G$, its \emph{boxicity} $\boxi(G)$ is the minimum dimension $b$, such that
$G$ is representable as the intersection graph of 
(axis--parallel) boxes in $b$--dimensional space. It is easy
to see that the class of graphs with $b\le 1$ is exactly the class of
interval graphs.
A $b$--dimensional \emph{box representation} of a graph $G=(V,E)$ is
a mapping  that maps each
each $u \in V$ to an axis-parallel $b$--dimensional box $B_u$
such that $G$ is the intersection graph of the family $\{B_u : u \in V\}$.

The concept of boxicity was  introduced by F. S.  Roberts \cite
{Roberts}. 
It  finds applications in niche overlap (competition) in ecology and to
problems of fleet maintenance in operations research. (See \cite {CozRob}.) 
It was shown 
by Cozzens \cite {Coz}  that computing the boxicity of a graph is NP--hard. 
This
was later strengthened by Yannakakis \cite {Yan1}, 
and finally by  Kratochvil \cite {Kratochvil} 
who  showed that deciding  whether
boxicity of a graph is at most 2 itself 
is NP--complete. The complexity of finding
the maximum independent set in bounded boxicity graphs was considered
by \cite {ImaiAsano,Fowler}.

There have  been many attempts to estimate or bound the boxicity of 
graph classes  with special structure. In his pioneering work,
F. S. Roberts proved that the boxicity of complete $k$--partite graphs
are $k$.  Scheinerman \cite {Scheiner} 
showed that the boxicity of outer planar graphs is at most $2$.
Thomassen \cite {Thoma1} proved that 
the boxicity of planar graphs is
bounded above by $3$. The boxicity of split graphs is investigated by
Cozzens and Roberts \cite{CozRob}.
Apart from these results, not much is known about
the boxicity of most of the well-known graph classes.
Also, little is known about the structure imposed on a graph by its high boxicity.

Researchers have also tried to generalize or extend  the
concept of boxicity in various ways. The poset boxicity \cite {TroWest}, 
the rectangular number \cite {ChangWest}, grid dimension \cite {Bellantoni},
circular dimension \cite {Feinberg,Shearer}  and the boxicity 
of digraphs \cite {ChangWest1} are some  examples. 

\subsection{Treewidth}
The notions of \emph{tree--decomposition} and \emph{treewidth}  
were first introduced (under different names) by R. Halin and later
rediscovered independently by Robertson and Seymour. (See
\cite {Diest}, Chapter 12  for historical details.) 
Roughly speaking, the 
treewidth of a graph $G$ is the minimum $k$ such that $G$ can be 
decomposed into pieces forming a tree  structure with at most $k+1$ 
nodes per piece. Such a decomposition is called a tree  
decomposition.  See section \ref {sec:defs} for the formal definition
of tree decomposition and treewidth.

These notions underly several important and  sometimes very deep results in 
graph theory and graph algorithms and are very useful for the analysis of 
several practical problems. 
Recent research has shown that many NP--complete problems become
polynomial or even linear time solvable, or belong to NC, when restricted
to graphs with small treewidth (See \cite {Arn,ArnPro,Bodland3}). 
The concepts of
treewidth and pathwidth have applications in many practically important
fields like VLSI layouts,  Cholesky factorization,
Expert systems,  Evolution theory, 
and natural language 
processing. (See  \cite {Bodland3} for references).

The decision problem of checking whether $\tw(G)$ is at most $k$, given
$G$ and $k$ 
is known to be NP-complete. Hence the   
problem of determining the treewidth of an arbitrary graph is NP-hard 
and the research on determining the  
treewidth and pathwidth has been focused on special classes.  
Linear or polynomial time or NC algorithms for producing optimal 
tree decompositions have been proposed for several 
special classes of graphs like graphs of bounded treewidth 
\cite{Bodland2,BodHag}, 
chordal graphs, cographs,  circular arc graphs, 
chordal bipartite graphs, permutation graphs, circle
graphs, and distance hereditary graphs.  
For an extensive bibliography on treewidth, see \cite {Bodland3}.

\section{Our Results}

Our main result is the following theorem which connects the boxicity of a graph
to its treewidth.

\noindent
{\bf {Theorem.}}
 \emph{ For any graph $G$, 
$\boxi(G) \, \le \, \tw(G) + 2.$
Moreover, we construct a family of graphs such that, for any $t \ge 1$, 
there is a graph $G$ of treewidth
at most $t + \sqrt{t}$ and boxicity at least $t - \sqrt{t}$ in this family.
In other words, for this family of graphs, 
$ \tw(G) (1 - o(1)) \le \boxi(G) \le \tw(G) + 2,$
and thus the upper bound is almost sharp.
}

As far as we know, the only known general upper bound for boxicity of $G$ is given
by Roberts\cite{Roberts}, who showed that $\boxi(G) \le \frac{n}{2}$, where $n$ is the
number of vertices in $G$. It was also shown that this bound is tight. Thus it is 
interesting to look for  upper bounds on boxicity which can
provide better structural insight about graphs with respect to their boxicity.
We believe that our upper bound in terms of treewidth is a progress in that direction. 
Since treewidth is an extensively studied concept, our upper bound
also leads to many nonintuitive results about boxicity, which appear difficult
to prove using  direct approaches.

\noindent
{\bf {Consequences on special graphs.}}
Except for the few classes mentioned in the introduction, not much progress is made
in bounding the boxicity of various well-known graph classes, perhaps due
to the lack of effective approaches. As consequences of the above theorem,
we are able to derive tight (up to constant factors) upper bounds for the 
boxicity of various graph classes, in terms of their maximum degree and clique number.

For a graph $G$, let $\Delta(G)$ denote its maximum degree and let $\omega(G)$
denote its clique number i.e., the size (number of nodes) of the maximum clique  in $G$.
We summarize our results for the boxicity of different graph classes in the following table.

\begin{tabular}{|c|c|}
\hline
Graph class & Upper bounds on $\boxi(G)$\\ \hline \hline
Chordal Graphs & $\omega(G) + 1$\\ 
  & $\Delta(G) + 2$\\ \hline
Circular Arc Graphs & $2\omega(G) + 1$ \\
 & $2\Delta(G) + 3$ \\ \hline
AT-free Graphs & $3\Delta(G)$ \\ \hline
Co-comparability Graphs & $2 \Delta(G) + 1$ \\ \hline
Permutation Graphs & $2 \Delta(G) + 1$ \\ \hline
Any minor closed family   & constant \\ 
 which excludes at least one planar graph  &  \\ \hline
\end{tabular}

\medskip
Each of the above upper bounds is shown to be tight up to small constant factors.
\vspace*{0.2cm}

\noindent
{\bf {Planar Graph Minors and Boxicity.}}
Study of graph minors is one of the most important areas in modern graph
theory (see \cite{Diest} for the definition of graph minors). 
Combing our upper bound result with a result of Robertson and Seymour \cite{RS13},
we obtain the following.

\noindent
{\bf {Theorem.}}
For every planar graph $H$, there is a constant $c(H)$ such that every graph
with boxicity $\ge c(H)$ has a minor isomorphic to $H$.

\noindent
{\bf {Cycles and Boxicity.}}
The properties that imposes various kinds of long cycles in a graph
is an extensively explored topic in graph theory.
(See Chapter 1 of \cite{Bondy_chap} for
an introductory survey, or the book by Voss \cite{Voss}.)
A structural consequence of our main result is that, high boxicity imposes
a long simple cycle as well as a long induced cycle (chordless cycle) in the graph.

\noindent
{\bf {Theorem.}}
\emph{In any  graph $G$ of boxicity $b$, there exists a simple cycle of length at least $b - 3$.
Moreover, there exists a graph $G$ whose boxicity is $b$ but the length of any simple cycle in it 
is at most $2b$. }

\noindent
{\bf {Theorem.}}
\emph{Let $G$ be a graph with maximum degree $\Delta$ and boxicity $b \ge 3$.
Then there exists a induced cycle (chordless cycle) of length at least
$\lfloor \log_{\Delta}(b - 2) \rfloor + 2$.}

\noindent
{\bf {Boxicity and Vertex Cover.}}
The subset $S \subseteq V(G)$ is called a \emph{vertex cover} of $G$
if every edge of $G$ is incident on at least one vertex from $S$.
A vertex cover of minimum cardinality is called a \emph{minimum vertex
cover}. We denote the cardinality of a minimum vertex cover of $G$  by 
$M\!V\!C(G)$. Our upper bound theorem yields the following.

\noindent
{\bf {Theorem.}}
For any graph $G$, $\boxi(G) \le M\!V\!C(G) + 2$.

The reader may note that for a graph $G$ on $n$ nodes, $M\!V\!C(G) = n - \alpha(G)$, where 
$\alpha(G)$ is the cardinality of a maximum independent set in $G$. Thus 
$\boxi(G) \le n - \alpha(G) + 2$. Similarly, $2$ times the cardinality of any maximal matching
in $G$ is an upper bound for $M\!V\!C(G)$ and thus, we obtain an upper bound for the
boxicity of $G$ in terms of the cardinality of maximal matchings also. 

We also connect boxicity to a parameter, which is a 
variant of minimum vertex cover, namely the cardinality of a minimum feedback vertex cover of $G$.

\noindent
{\bf {Algorithmic Consequences.} }
From an algorithmic point of view, it is interesting to efficiently construct 
a box representation (see introduction for the definition of box representation) 
of the given graph in low dimensional space. For example,
note that if a dense graph has such a representation in constant dimensional space,
then the memory required to store this graph (nodes and edges) is only linear 
in the number of nodes.
The proof of our main upper bound also yields an efficient construction
of the box representation of the given graph $G$ from its tree decomposition.
Since, for bounded treewidth graphs, the tree decompositions can be constructed
in linear time \cite{Bodland2}, we have the following.

\noindent
{\bf {Theorem.}}
\emph{For a  bounded treewidth graphs,
box representation in constant dimension can be constructed
in linear (in the number of vertices) time.}

Efficient polynomial time 
algorithms (exact or approximation) are known for constructing
the tree decompositions of many special graph classes such as 
chordal graphs, cographs, circular arc graphs, chordal bipartite graphs,
permutation graphs, circle graphs, distance hereditary graphs etc.
An immediate consequence  of this in conjunction with the constructive 
proof of our main upper bound
 is that the corresponding box representations
can also be computed in polynomial time.

\medskip

\noindent
\textbf{Complexity theoretic consequences:} By our main upper bound result, 
the class of bounded treewidth graphs is a subset of bounded boxicity graphs.
Hence if a problem is NP-hard for bounded treewidth graphs, it is also NP-hard
for bounded boxicity graphs. For example, it is shown in \cite{McDiarmid} that
the channel assignment problem is NP-complete for graphs of treewidth at least $3$.
It follows from our result that this problem is NP-complete for graphs of boxicity
at least $5$.

\section{Tree Decompositions and the Treewidth}
\label{sec:defs}

\begin {defn}\label{def:tree-decomp}
A \textbf{tree decomposition} of $G=(V,E)$ is  a pair 
$(\{X_i : i \in I\}, T),$  where  $I$ is an index set, $\{X_i: i \in I\}$ is a collection of
subsets of $V$ 
and $T$ is a tree (connected) whose node set is $I$, 
such that the following conditions are satisfied:
\begin {enumerate}
\item  $\bigcup_{i \in I} X_i = V$.
\item  $ \forall (u,v) \in E, \exists i \in I \mbox{~~such that~~}  u,v \in X_i$.
\item  $ \forall i,j,k \in I$: if $j$ is on a path in $T$
       from $i$ to $k$, then $X_i \cap X_k \subseteq  X_j$.
\ignore{
 This property can be equivalently stated as follows: 
For a node  $x \in V$,  let $I_x \subseteq I$ be the set of nodes of $T$
such that $x \in X_i$ if and only if $i \in I_x$.
, where 
 $I_x \subseteq I$ is the set of nodes of $T$ such that $x \in X_i$ if and
 only if $i \in I_x$. }
\end {enumerate}

\noindent
The \textbf{width} of a tree decomposition $(\{X_i: i \in I \}, T)$
is  $\max_{i \in I} |X_i| - 1$.  
The \textbf{treewidth} of $G$ is the minimum width over all
tree decompositions of $G$ and is denoted by $\tw(G)$. 
Node $i$ of a tree decomposition $(\{X_i: i \in I\}, T)$ refers to the node $i$ of the tree $T$.
\ignore{
When there is no ambiguity, when we say node $X_i$ of the tree decomposition,
it refers to the 
nodes $X_i$ as  node $i$ of the 
We also call the subset of nodes $X_i$ as the node $i$ of the tree
decomposition when it is clear from the context.}

\end {defn}

\noindent
\emph{\textbf{Rooted Tree.}} A tree with a fixed root is called a rooted tree.
The  $height(i)$ of a node $i$ in a rooted tree $T$ with root $r$ 
is defined as usual:  $height(r)$ of the root $r$ is 0, 
and $height(x)$ for any other node $x$ is  exactly one more than the
height of its parent. 
A node $i \not= j$ is the \textbf{ancestor} of node
$j$ if $i$ is in the path from $j$ to $r$.  A node $j$ is a \textbf{descendant}
of $i$ if either $i=j$ or $i$ is the ancestor of $j$.

\begin {defn} \label{def:norm-tree-decomp}
A \textbf{normalized tree decomposition} of a graph $G = (V, E)$ is a triple $(\{X_i:
i \in I\},\, r \in I,\, T)$ where $(\{X_i : i \in I\},\, T)$ is  a tree decomposition of $G$  that
additionally satisfies
the following two  properties. 
\begin{enumerate}
\item[4.]
  It is a rooted tree where the subset $X_r$ that corresponds  to
  the root node $r$ contains exactly one vertex. 

\item[5.]
  For any node $i$, if $i'$ is a child of  $i$, then  $|X_i' - X_i| = 1$.

\end{enumerate}

\end {defn} 

\begin{lem} \label{lemma:norm-tree}
  For any graph $G$ there is a normalized tree decomposition with width  equal
  to  $\tw(G)$.
\end{lem}

\begin{pf}
  Consider a tree decomposition $(\{X_i: i \in I\}, T)$ of  $G =(V,E)$ with width $\tw(G)$.
  We  convert  it into a
  normalized 
  tree decomposition $(\{X_i: i \in I'\},\, r,\, T')$ as follows. 

  As the first step, we convert $T$ into a rooted tree $T_1$ as follows.
  Let $i$ be an arbitrary node of $T$ such that $X_i$ is non--empty. Let 
  $u \in X_i$. Create a new node $r$ (where $r \notin I$), and define $X_r = \{u\}$. Now connect
  node $r$ to  $i$. Let the resulting tree on the node set
  $I \cup \{r\}$ be $T_1$. It is easy to verify that $(\{ X_i : i \in I \cup
  \{ r\}\},\, T_1)$  is a tree decomposition of $G$.
  From here on, we view $T_1$ as a rooted tree, with root $r$. 
  
  Consider any edge $(j,j')$ of $T_1$  where
  $j'$  is  a child of  $j$. Without loss of generality, 
  we can assume that $X_{j'} \not \subseteq
  X_j$. (If $X_{j'} \subseteq X_j$ then the following operations do not violate
  the defining properties of tree decomposition: (a) Remove $j'$ from  $I
  \cup \{r\}$ and hence from $T_1$\, (b) make each child of $j'$ a child of
  $j$.)
  Let $X_{j'} - X_{j}  = \{u_1,\ldots,u_h\}$.
  If $h  = 1$  then we retain this edge as such.
  If $h > 1$ then we replace the edge $(j,j')$ by a path
  $j, k_1, k_2, \cdots, k_{h-1},j'$, where
  $k_1,k_2, \cdots, k_{h-1}$, are new nodes, and define the subset 
  $X_{k_i} =(X_j \cap X_{j'}) \cup  \{u_1, u_2, \cdots,u_i\}$  for $1 \le i \le h-1$. 
  Note that $|X_{k_i}| \le |X_{j'}|$ for $1 \le i \le h-1$ 
  and thus by introducing these new nodes we have not increased the width of
  the tree decomposition. We repeat this process for each edge of $T_1$. 
  Let $T'$  be the new rooted tree (rooted at $r$)  
  obtained after these operations.  Let $I'$ be the node set of $T'$.
  Note that the root $r$ of $T'$ still corresponds to
  the singleton set $X_r$.
  Now, it is straightforward  to verify that $(\{X_i: i \in I'\},\, r,\, T')$
  is a normalized tree decomposition. 
\end{pf}

\medskip
\noindent
\begin {defn}\label{def:bh}
With respect to the normalized tree decomposition $\nr{X}{I}{r}{T}$  of a graph $G=(V,E)$,
we define the following two functions.

\begin {enumerate} [a)]

\item  $\bag : V \rightarrow I$  is defined as follows. For $v \in V$,  $\bag(v) = i$, where
 $i$ is the  (unique) node in $I$ such that $height(i)$ is minimum subject to
 the condition that  $v \in X_i$.

\item $\level: V \rightarrow N$ is defined by $\level(v) = height(\bag(v))$. 

\end {enumerate} 
\end{defn}

Observe that the function $\bag(v)$ is well-defined. That is, there is exactly one
node $i$  of $T$, such that  $v \in X_i$ and $height(i)$ is the minimum
possible. To see this, assume that there is one more node $j$ such that $v \in X_j$ and
$height(j) = height(i)$. Then, by 
Property 3 of Definition \ref{def:tree-decomp},  there should be a node $k$
with height less than $height(i)$ and $v \in X_k$. This
contradicts the assumption that node $i = \bag(v)$ has the minimum possible
height.

\begin {lem} \label{lemma:bijective}
  	The function $\bag : V \rightarrow I$  is 
	a bijection. 
\end {lem} 
  \begin {pf}
  First  
  we show that $\bag : V \rightarrow I$ is injective. 
  That is, for any two distinct vertices  $u,v \in V$,~ $\bag(u) \not=  \bag(v)$.
  If not, let  $\bag(u) = \bag(v) = i$. Since $X_i$ contains at least
  two vertices (namely $u$ and $v$),  $i$ is not the root node of the normalized tree
  decomposition.
  Let $j$ be the parent of $i$.  Since
  $\bag(u) = \bag(v) = i$, $X_j$ does not contain $u$ and $v$ by the definition of
 $\bag(\cdot)$. That is, $|X_i - X_j| \ge 2$.
  This  contradicts Property 5  of Definition \ref{def:norm-tree-decomp}.
  
  Now assume that $\bag : V \rightarrow I$
  is not surjective i.e., there exists a node  $i \in I$ that do not have a
  pre-image in $V$.  Let $X_i = \{u_1,\ldots,u_r\}$.
  Consider any vertex $u_j \in X_i$. Let $\bag(u_j) = k \not= i$. By the
  definition of $\bag(u_j)$, $height(k) \le height(i)$. Let $j$ be the parent
  of $i$. Clearly, $j$ is on the path  between $k$ and $i$ in $T$, and thus
  by Property 3 of the Definition \ref{def:tree-decomp}, $u_j \in X_j$.
  This implies that $X_i \subseteq X_j$, which
  contradicts Property 5 of Definition \ref{def:norm-tree-decomp}.
  Thus $\bag : V \rightarrow I$ injective and surjective i.e., bijective.
  \end {pf}

\begin{lem}\label{lemma:decent}
  For any $i \in I$ such that $u \in X_i,$ node $i$ is a descendant of $\bag(u)$.
\end{lem}
\begin{pf}
  Otherwise, since $u \in X_i$ as well as $X_{\bag(u)}$,
  by Property 3 of Definition \ref{def:tree-decomp},  $u \in X_j$ also where
  $j$ is the parent of $\bag(u)$. This contradicts the definition of
  $\bag(u)$. 
\end{pf}


\section{Box Representation and  Interval Graph Representation}
\label{sec:boxi-inter}

Let $G=(V,E(G))$ be a graph and let $I_1, \ldots, I_k$ be $k$
interval graphs such that each $I_j = (V, E(I_j))$ is defined on
the same set of vertices $V$.
If $$E(G) = E(I_1) \cap \cdots \cap E(I_k),$$ then we say that
$I_1, \ldots, I_k$ is an \emph{interval graph representation} of $G$.
The following equivalence is well-known.

\begin{thm}[Roberts \cite{Roberts}]\label{thm:robertsx}
The minimum $k$ such that there exists an interval graph representation
of $G$ using $k$ interval graphs $I_1, \ldots, I_k$  is the same
as $\boxi(G)$. 
\end{thm}

Recall that a $b$--dimensional box representation of $G$ is a mapping of each node $u \in V$ 
to $R_1(u) \times \cdots \times R_b(u)$, where each $R_i(u)$ is a closed interval of the
form $[\ell_i(u), r_i(u)]$ on the real line. It is straightforward to see that 
an interval graph representation of $G$ using $b$ interval graphs $I_1, \ldots, I_b$,
is equivalent to a $b$--dimensional box representation in the following sense. 
Let $R_i(u) = [\ell_i(u), r_i(u)]$ denote the closed interval corresponding to node $u$
in an interval realization of $I_i$. 
Then the $b$--dimensional box corresponding  to $u$  is simply
$R_1(u) \times \cdots \times R_b(u)$. Conversely, given a $b$--dimensional box representation
of $G$, the set of intervals $\{R_i(u) : u \in V\}$ forms the $i$th interval graph $I_i$
in the corresponding interval graph representation.

\newcommand \edgeprop {Edge Property}
\newcommand \connprop {Connectivity Property}

\ignore{
	\section{Treewidth vs Boxicity}
	
	The proofs of the following theorems are involved and requires some
	conceptual development. Due to lack of space, these proofs are moved
	to the Appendix.
	
	\begin{thm}\label{thm:tw-upper}
	  For any graph $G$, $\boxi(G) \, \le \, \tw(G) + 2$.
	\end{thm}
	
	When we consider various simple examples, it is tempting to conjecture
	that the tight upper bound on the boxicity is $\frac{\tw(G)}{2}$.
	(For example, consider the Roberts graph explained in Section \ref{sec:circarc}).
	But we show that the above upper bound is asymptotically tight.
	More precisely,
	
	\begin{thm}\label{thm:tw-lower}
	For any integer $k \ge 1$, there exists a graph $G$ with $\tw(G) = k$ and
	$\boxi(G) \ge k (1 - \frac{2}{\sqrt{k}}) = k (1 - o(1))$.
	\end{thm}
	
	The proof of Theorem \ref{thm:tw-upper} is constructive.
	It is straight forward to verify from the proof that,
	 starting with a tree decomposition 
	of constant width, the box representation of the graph can be constructed in linear
	time. Bodlaender \cite{Bodland2} proved that for bounded treewidth graphs, the optimum tree decomposition
	can be constructed in linear time.

	\begin{thm}
	For a  bounded treewidth graph $G$ on $n$ nodes,
	a  box representation of $G$ in constant dimension can be constructed
	in $O(n)$ time.
	\end{thm}

}

\section{Treewidth vs Boxicity: The Upper Bound} \label{sec:up}

Let $G=(V,E)$ be a graph. 
In this section, we assume that $\nr{X}{I}{r}{T}$ is a normalized tree decomposition 
of $G$ with width $\tw(G)$.

\begin {lem} \label{lemma:theta}
  Let $G=(V,E)$ be a graph and $\nr{X}{I}{r}{T}$ be its normalized tree decomposition
  of width $\tw(G)$. 
  Then, there exists a function $\theta : V \rightarrow \{0, \ldots,
  \tw(G) \}$,
  such that for any $i \in I$ and for any two distinct nodes  $u, v \in X_i,$ $\theta(u) \not=
  \theta(v)$.
\end {lem}

  \begin{pf}
  Sort the nodes in $I$ in the increasing order of their height
  (breaking ties arbitrarily). Let the order be $i_1,\ldots,i_n$.
  Let $u_j = \bag^{-1}(i_j)$. 
  We inductively define $\theta(u_j)$ in
  the order $u_1, \ldots, u_n$. 
  Define $\theta(u_1) = 0$.
  Assume inductively that for $k < j$,
  $\theta(u_k)$ is defined, and for any $u, v \in X_{{i_k}}$, $\theta(u)
  \not= \theta(v)$. 
  Observe that node $i_1$ is the root of $T$ and
  thus $X_{i_1}$ is the singleton set $\{u_1\}$, and thus the inductive
  assumption is trivially true for $i_1$. Let $i_h$ be the parent of $i_j$ in $T$.
  First we observe that $u_j \in X_{i_j} - X_{i_h}$, by the
  definition of $\bag(u_j)$. 
  Hence $X_{i_j} - \{u_j\} \subseteq X_{i_h}$ by Property 5 of Definition
  \ref{def:norm-tree-decomp}.
  Consider a vertex $v \in  X_{i_j} - \{u_j\}$.
  Observe that $\bag(v) = i_r$, for some  $r < j$.
  (This is because, $v \in X_{i_1} - \{u_j\} \subseteq X_{i_h}$ and $height(i_h) < height(i_j)$.)
  So $\theta(v)$ is already defined at this point by the induction assumption. Now define $\theta(u_j) =
  t$, where $t \not= \theta(v)$ for any $v \in X_{i_j} - \{u_j\}$.
  There is such a $t$ because $|X_{i_j} - \{u_j\}| \le \tw(G)$ but there
  exists $\tw(G) + 1$ distinct possible values for $t$.
  Now we claim that for any $u, v \in X_{i_j} - \{u_j\}$,  $\theta(u) \not=
  \theta(v)$. This is because, $u,v \in X_{i_j} - \{u_j\} \subseteq X_{i_{h}}$, and since $h < j$, 
  the inductive assumption is valid for $i_h$.
  \end{pf}

  \begin{lem}\label{lemma:grouping-1}
    If $(u, v) \in E(G)$ then $\theta(u) \not= \theta(v)$.
  \end{lem}
  \begin{pf}
  	Since $(u, v) \in E(G)$ then there exists an $X_i$ such that
	$u,v \in X_i$ by Property 2 of Definition \ref{def:tree-decomp}.
	Now, by Lemma \ref{lemma:theta}, $\theta(u) \not= \theta(v)$.
  \end{pf}

\begin{lem}\label{lemma:prop-norm-tree}
  If  $(u,v) \in E(G)$,  then either $\bag(u)$ is an ancestor of $\bag(v)$ or
  $\bag(v)$ is an ancestor of $\bag(u)$ in $T$.
\end{lem}

\begin{pf}
  Due to Property 2 of Definition \ref{def:tree-decomp},  there is a node
  $i \in I$ such  that $u, v \in X_i$. Because of Lemma \ref{lemma:decent},
  node $i$ is the descendant of $\bag(u)$ and also $\bag(v)$.
  Thus both $\bag(u)$ and $\bag(v)$ are in the path from $i$ to the root $r$ in $T$.
  Moreover, $\bag(u) \not= \bag(v)$ since $\bag(\cdot)$ is a bijection  by Lemma 
\ref{lemma:bijective}.
  Thus the result follows.
\end{pf}

  \begin{lem}\label{lemma:grouping}
  	Let $(u,v) \in E(G)$ and let $\bag(u)$ be the ancestor of \,$\bag(v)$. 
	For any vertex $w \in V - \{u\}$,
	$\theta(w) \not= \theta(u)$ if $\bag(w)$ is in the path
	from $\bag(v)$ to $\bag(u)$ in $T$.
  \end{lem}
  \begin{pf}
	Because of Property 2 of Definition \ref{def:tree-decomp}, there is an
	$X_i$ such that $u, v \in X_i$. By Lemma \ref{lemma:decent},
	we know that $i$ is a descendant of $\bag(u)$ and also $\bag(v)$.
	This in conjunction with the assumption that $\bag(u)$ is the ancestor of $\bag(v)$, 
	implies that $\bag(v)$ is in the 
	path from $\bag(u)$ to node $i$. Thus, for any node $k$ in the path
	from $\bag(v)$ to $\bag(u)$, $u \in X_k$, by Property 3 of 
	Definition \ref{def:tree-decomp}. Now, for any vertex $x \in X_k - \{u\}$,
	~$\theta(x) \not= \theta(u)$ by Lemma \ref{lemma:theta}. In particular,
	this is true for $k = \bag(w)$ and $x=w$.
  \end{pf}

\medskip
\noindent
Using the function $\theta : V \rightarrow \{0, \ldots, \tw(G)\}$ (see Lemma  \ref{lemma:theta}) and function $h : V \rightarrow N$ (see Definition \ref{def:bh}), we construct $tw(G) + 2$ different interval super
graphs of $G$ as follows.
Let $i$ be such that $0 \le i \le \tw(G)$. The interval graph $I_i$ is defined
as follows.

  \medskip
  \noindent {\bf Definition of interval graph $I_i,$~ for $0 \le i \le \tw(G)$:}
  We define the interval $[\ell_i(v), r_i(v)]$ for each $v \in V$ as follows.
  \begin {enumerate}
  
  \item  If $\theta(v) = i$ then $\ell_i(v) = 2 \level(v)$ and $r_i(v) = 2\level(v) +
  1$.
  
  \item  If $\theta(v) \not= i$ then let $S = \theta^{-1}(i) \cap N(v)$, where
  $\theta^{-1}(i) = \{u \in V ~|~ \theta(u) = i\}$ and $N(v) = \{u \in V -
  \{v\} ~|~ (u, v) \mbox{~is an edge of~} G\}$.
  \begin {enumerate} 
   \item If $S = \emptyset$ then $\ell_i(v) = 3n$ and $r_i(v) = 3n$.
    \item If $S \not= \emptyset$ then $\ell_i(v) = \min_{u \in S} r_i(u)$ and
	$r_i(v) = 3n$.
  \end {enumerate}

\end {enumerate}

  \medskip
  \noindent {\bf Definition of interval graph  $I_{\tw(G) + 1}$:}
  Consider a \emph{depth-first} ordering of the nodes of $T$. The depth-first
  ordering of rooted tree $T$ rooted at $r$ is an ordered list of the nodes of $T$
  denoted as $df(T, r)$. If $T$ has only one node,
  namely its root $r$, then $df(T,r) = \langle r, r \rangle$.  Otherwise, let $r_1, \ldots,
  r_k$ be the children of $r$ and let $T_i$ be the rooted sub-tree rooted at
  $r_i$. Then, $df(T,r)$ is the concatenation of the lists $\langle r \rangle, df(T_1,r_1), \ldots, df(T_k, r_k), \langle r
  \rangle$ in that order.
  Observe that each node of $T$ appears exactly two times in $df(T, r)$.
  Thus we can associate with each node $i$, two numbers $first(i)$ and $last(i)$
  that denote its sequence number in the ordered list $df(T, r)$ corresponding to its first
  occurrence and last occurrence respectively.
  Now, for each vertex $v \in V$, $\ell_{\tw(G)+1}(v) = first(\bag(v))$ and $r_{\tw(G)+1}(v) =
  last(\bag(v))$.  The resulting interval graph is $I_{\tw(G)+1}$.

  \begin{lem} \label{lemma:edge}
  Each $I_i ,$ ~$0 \le i \le \tw(G) + 1$,  $E(G) \subseteq E(I_i)$.
  \end{lem}

 \begin{pf}
  Let  $(x, y) \in E(G)$. First, assume that $0 \le i \le \tw(G)$.
  By Lemma \ref{lemma:grouping-1}, we have  $\theta(x) \not= \theta(y)$.
  Without loss of generality, assume that $\theta(y) \not= i$. Hence 
  $r_i(y) = 3n$ because of Case 2 of the definition of $I_i$.
  If $\theta(x) \not= i$ then $r_i(x)$ is also $3n$, and thus $(x, y) \in E(I_i)$.
  Now assume that $\theta(x) = i$.  Hence,
  $x \in S = \theta^{-1}(i) \cap N(y)$. Hence 
  $
  	r_i(y) = 3n \ge r_i(x) \ge \min_{z \in S} r_i(z) = \ell_i(y),
  $
  and thus $r_i(x) \in [\ell_i(y), r_i(y)]$. It follows that $(x, y) \in E(I_i)$.

    It remains to show that $(x, y) \in E(I_{\tw(G) + 1})$.
	Because of Lemma \ref{lemma:prop-norm-tree},
	we can assume without loss of generality that $\bag(x)$ is an ancestor of $\bag(y)$.
	Consider the depth-first order $df(T, r)$ of the nodes of $T$.
	It is straightforward to verify that $first(\bag(x)) \le first(\bag(y))
	\le last(\bag(y)) \le last(\bag(x))$, and thus the intervals in $I_{\tw(G)
	+ 1}$ corresponding
	to $x$ and $y$ intersect. It follows that $(x, y) \in E(I_{\tw(G) + 1})$.
   \end{pf}

   \begin{lem} \label{lemma:no-edge}
   For any  $(x, y) \notin E(G)$, there exists some $i, ~0 \le i
   \le \tw(G) + 1$, such that $(x, y) \notin E(I_i)$.
   \end{lem}
	 
  \begin{pf}
   Let $(x, y) \notin E(G)$. Suppose that neither $\bag(x)$ is an ancestor of $\bag(y)$
   nor $\bag(y)$ an ancestor of $\bag(x)$ in $T$. Then, we claim that $(x, y) \notin
   E(I_{\tw(G)+1})$. This is because, if $\bag(x)$ is not an ancestor of
   $\bag(y)$ or vice versa, then in the depth-first
   order $df(T, r)$, either $last(\bag(x)) < first(\bag(y))$ or
   $last(\bag(y)) < first(\bag(x))$, and thus their corresponding intervals do
   not intersect.
   \emph{From now on, we assume without loss of generality, that $\bag(x)$ is an ancestor of
   $\bag(y)$.}  
   
   Let $t = \theta(x)$. Note that by the definition of $\theta$,~ $0 \le t \le
   \tw(G)$. We claim that $(x, y) \notin E(I_t)$. Since function $\bag : V
   \rightarrow I$ is bijective (Lemma \ref{lemma:bijective}), $\bag(x) \not=
   \bag(y)$, and thus, since $\bag(x)$ is an ancestor of $\bag(y)$, we have
   $\level(x) < \level(y)$. 
   
   Now, if $\theta(y)$  also equals $t$, 
   then the intervals corresponding to $x$ and $y$ do
   not intersect
   since $\level(x) \not= \level(y)$ (see definition of $I_t$).
   \emph{From now on, we assume that $\theta(y) \not= \theta(x)$.}

   Let $S = \theta^{-1}(t) \cap N(y)$. If $S = \emptyset$ then the interval
   corresponding to $y$ in $I_t$ is $[3n, 3n]$ by definition. Since $r_t(x) = 2 \level(x) +1 < 3n,$~~
   $(x, y) \notin E(I_t)$ as required. If $S \not= \emptyset$ then let
   $z$ be the node such that $r_t(z) = \min_{w \in S} r_t(w)$. 
   Note that $x \not= z$ since $x \notin N(y)$. Since $z \in N(y)$,
   there is an ancestorial relation between $\bag(z)$ and $\bag(y)$ by Lemma
   \ref{lemma:prop-norm-tree} i.e., $\bag(z)$ is an ancestor of $\bag(y)$ or vice versa.
   Recalling that $\bag(x)$ is an ancestor of $\bag(y)$, it follows that 
   \emph{there is a pair-wise ancestorial relation between $\bag(x)$, $\bag(y)$
   and $\bag(z)$}. Noting that $x \not= y \not= z$, by Lemma \ref{lemma:bijective} $\bag(x) \not= \bag(y) \not=
   \bag(z)$. It follows that $\level(x) \not= \level(y) \not= \level(z)$.
   
   Let $\level(x) > \level(z)$.
   Recalling that $\level(y) > \level(x)$, we have
   $\level(y) > \level(x) > \level(z)$
   Hence $\bag(x)$ is in the path in $T$ from $\bag(y)$ to $\bag(z)$, and thus
   by Lemma \ref{lemma:grouping}, it follows that $\theta(x) \not= \theta(z)$. 
  But recall that, $z \in S \subseteq \theta^{-1}(t)$, by definition, and
 therefore $\theta(z) = t = \theta(x)$, which is a contradiction.
Hence the only possibility is that 
   $\level(x) < \level(z)$, and thus $r_t(x) < r_t(z)$ by case (1) of
   definition of $I_t$.
   Recall that  $r_t(z) = \ell_t(y)$ by case 2(b) of definition of $I_t$. 
	Hence we have $r_t(x) < \ell_t(y)$,
   and thus
   $(x, y) \notin E(I_t)$.
  \end{pf}

\medskip
\noindent
By combining  Lemma \ref{lemma:edge} and Lemma
\ref{lemma:no-edge}, we can infer that $E(G) = E(I_1) \cap \cdots \cap E(I_{\tw(G) + 1})$.
Thus, by Theorem \ref{thm:robertsx}, we obtain the following.

\begin{thm}
\label{thm:tw-upper}
For any graph $G$, $\boxi(G) \, \le \, \tw(G) + 2$.
\end{thm}

\ignore{
	\begin{pf}
	
	  Let $\tw(G) = t$.
	  By Lemma xx, it is enough to show that there exists $t + 2$
	  interval graphs $I_1,\ldots,I_{t+2}$ interval graphs such that
	  $$
	  	G = I_1 \cap I_2 \cap \ldots \cap I_{t+2}.
	  $$
	  Consider a normalized tree decomposition $T$ of $G$. 
	  \medskip
	
	  Let $Q: V \rightarrow \{Q_1,\cdots, Q_{t+1}\}$ be a partitioning of the 
	  vertices of $G$, such that for each $i \in I$, $|X_i \cap Q_j| = 1$. 
	  (The existence of such a partitioning is ensured by Lemma \ref {}.) 
	
	  Now we construct $t+2$ interval graphs denoted as $I_1,\ldots,I_{t+1}$, and
	  $I_s$, where the first $t+1$ interval graphs are constructed based on the above
	  grouping, in a similar fashion, and the last `special' interval graph $I_s$
	  is constructed in a different manner. Later we will show how these $t+1$
	  interval graphs together realize~$G$.
	  \end{pf}

}


\section{Tightness Result}

When we consider various simple examples, it is tempting to conjecture
that the tight upper bound on the boxicity is $\frac{\tw(G)}{2}$.
(For example, consider the Roberts graph explained in Section \ref{sec:circarc}).
But we show that the above upper bound is asymptotically tight.
More precisely,

%

\begin{thm}\label{thm:tw-lower}
For any integer $k \ge 1$, there exists a graph $G$ with $\tw(G) \le k$ and
$\boxi(G) \ge k (1 - \frac{2}{\sqrt{k}}) = k (1 - o(1))$.
\end{thm}

\begin{pf}

We show the following. Fix any $t \ge 1$. We construct a graph $G$ such that
$\tw(G) \le t + \sqrt{t}$ ~and ~$\boxi(G) ~\ge~ t - \sqrt{t}$.
For any fixed $k$, we get the result by choosing $t$ to be the largest
integer such that $t + \sqrt{t} \le k$.

The graph $G$ is as follows.
The node set of $G$ is the disjoint union of 
 $\alpha + 1$ sets  $P_0, P_1,
\ldots, P_\alpha,$ where 
$$
\alpha = \sum_{j=1}^{\lfloor\sqrt{t}\rfloor} {t \choose j}.
$$
(Note that $\alpha$ corresponds to the total number of 
non-empty subsets of at most $\sqrt{t}$ elements from a collection of
$t$ distinct elements.)
Let $|P_0| = t$. The cardinality of $P_i$ for $1 \le i \le \alpha$ is
defined as follows. 
Let $\mathcal{S} = \{ A \subseteq P_0 ~|~ 1 \le |A| \le \sqrt{t}\}$.
Note that $|\mathcal{S}| = \alpha$.
Let $\Pi: \mathcal{S} \rightarrow \{P_1, \ldots, P_\alpha\}$ be a bijective 
map.  We define $|\Pi(A)| = |A|$. In other words, $|P_i| = |\Pi^{-1}(P_i)|$.

For $i \in \{1, \ldots, \alpha\}$, let $\core_i: P_i \rightarrow \Pi^{-1}(P_i)$
be a bijection.

Edge $(u, v) \in E(G)$ if and only if any one of the following conditions holds
\begin{enumerate}[a)]
\item
$u, v \in P_i$ for some $i$.
\item
$u \in P_0$ and $v \in P_i$ for some $i \ge 1$ and $u \not= \core_i(v)$.
\end{enumerate}

\begin{claim}
$\tw(G) \le t + \lfloor \sqrt{t} \rfloor$. 
\end{claim}
The above claim can be seen as follows.
Define a tree decomposition $(\{X_i : i \in I\}, T)$ of $G$ where
$I = \{0, \ldots, \alpha\}$. Define $X_0 = P_0$ and $X_i = P_0 \cup P_i$
for $i \in \{ 1, \dots, \alpha\}$.
The edge set of $T$ is $\{(0, i): 1 \le i \le \alpha\}$.

(Note that $T$ is a star with node $0$ at the center).
It is straightforward to verify that this is a valid
tree-decomposition of $G$. Recalling that each $P_i, 1 \le i \le \alpha,$ has
at most $\sqrt{t}$ nodes, it follows that the width of this decomposition is
$t + \floor{\sqrt{t}}$.

\begin{claim}
$\boxi(G) \ge t - \lfloor \sqrt{t} \rfloor$.
\end{claim}

Proof:
Assume by contradiction
that $\boxi(G) < t - \floor{\sqrt{t}}$. Then, consider an interval graph representation of
$G$ using $\gamma = t - \floor{\sqrt{t}} - 1$ interval graphs.
That is, let
$$E(G)= E(I_1) \cap \ldots \cap E(I_\gamma),$$
where $I_1, \ldots, I_\gamma$ are the $\gamma$ interval graphs.
Fix any arbitrary interval realizations for $I_1, \ldots, I_\gamma$.
From now on (abusing the terminology), we refer to this interval realization of $I_i$
also as $I_i$.

\medskip
\noindent

Let $\mathbf{P}: \{I_1, \ldots,
I_\gamma\} \rightarrow P_0 \times P_0$ be a function, where $\mathbf{P}(I_j)$ is
defined as follow.
For  node $w \in V$, let $[\ell_j(w), r_j(w)]$ denote its corresponding interval in $I_j$.
Let $u, v \in V$ be such that 
$\ell_j(u) = \max_{w \in P_0} \ell_j(w)$ and $r_j(v) = \min_{w \in P_0} r_j(w)$ (resolving ties
arbitrarily).
Define $\mathbf{P}(I_j) = (u, v)$.
Recalling that $P_0$ induces a complete graph in $G$, 
it follows that for any node $w \in P_0$, $r_j(w) \ge \ell_j(u)$ (otherwise intervals
corresponding to nodes $u$ and $w$ will not intersect). Thus $r_j(v) \ge \ell_j(u)$ 
and therefore $[\ell_j(u), r_j(v)]$ is a valid interval.
Now it is straightforward to see that 
for any node $w \in P_0,$ $[\ell_j(u), r_j(v)] \subseteq [\ell_j(w), r_j(w)]$.
Note that $[\ell_j(u), r_j(v)] = [\ell_j(u), r_j(u)] \cap [\ell_j(v), r_j(v)]$.
Now it is easy to see that there
cannot be a node $x \in P_k$, for any $k \ge 1$, that is adjacent to both $u$ and $v$
 but not to some $y\in P_0$, in
this interval graph. This is because, if $x$ is adjacent to both $u$ and $v$, then 
the interval corresponding to $x$
has a non-empty intersection with $[\ell_j(u), r_j(v)]$, and thus it has
a non empty intersection with the interval for any node $y \in P_0$.
This is summarized as follows.


\begin{claim}\label{claim:four-nodes}
Consider any interval graph $I \in \{I_1, \ldots, I_\gamma\}$.
Let $\mathbf{P}(I) = (u, v)$.
Let $x \in P_k$ for any $k \in \{1, \ldots, \alpha\}$.
If $\core_k(x) \notin \{u, v\}$ then edge $(x, \core_k(x)) \in E(I)$.
\end{claim}

Define a multi-graph $H=(V_H, E_H)$, where  
$$
V_H = P_0  ~~~~\mbox{and the multi-set}~~
E_H = \{\mathbf{P}(I_1), \mathbf{P}(I_2), \ldots, \mathbf{P}(I_\gamma)\}.
$$ 
Note that $H$ has $t$ nodes and $\gamma$ edges.

Applying  Lemma \ref{lemma:connected-component} to $H$ by fixing $r= \lfloor \sqrt{t} \rfloor$,
we infer the following.
The multi-graph $H$ has a connected component $K=(V_K, E_K)$ on $k$ nodes
and exactly $k-1$ edges,
where $1 \le k \le \floor{\sqrt{t}}$.
Let
$$
S = \{I ~~|~~ \mathbf{P}(I) \in E_K \}. \mbox{~~~~Clearly~} |S| = k-1.
$$



Let $\Pi(V_K) = P_k$. 
Define a function $\finter: P_k \rightarrow \{I_1, \ldots, I_\gamma\}$
such that for $x \in P_k$, $\finter(x) = I_j$ where $I_j$ is an interval graph
such that $(x, \core_k(x)) \notin E(I_j)$. (Reader may note that
there exists one such interval graph because $(x, \core_k(x)) \notin E(G)$. On the
other hand, there can be more than one interval graph where the edge $(x, \core_k(x))$ is not present.
The function $\finter(x)$ maps $x$ to one such interval graph.)

\begin{claim}\label{claim:finter}
For any $x \in P_k = \Pi(V_K),$~ $\finter(x) \in S$.
\end{claim}

Proof: Let $I \in \{I_1, \ldots, I_\gamma \} - S$ and let $\mathbf{P}(I) = (u, v)$.
Since $(u, v) \in E(H)$, both $u$ and $v$ belong to the same connected component of $H$.
We claim that $u, v \notin V_K$. Otherwise, $(u, v) = \mathbf{P}(I) \in E_K$, and
hence by definition of $S$, $I \in S$, a contradiction.
Recall that $P_k = \Pi(V_K)$ and thus  for any $x \in P_k$, $\core_k(x) \in V_K$ by definition
of the function $\core_k(\cdot)$.
Therefore $\core_k(x) \notin \{u, v\}$.
It follows by Claim \ref{claim:four-nodes}, that $(x, \core_k(x)) \in E(I)$ and therefore
$\finter(x) \not= I$ by the definition of $\finter(x)$. 
The claim follows.

Recall that $|V_K| = k$, where $1 \le k \le \floor{\sqrt{t}}$.
It follows that $|P_k| = k$ since $P_k = \Pi(V_K)$ and recalling that
$|\Pi(V_K)| = |V_K|$ by definition.
But recall that $|S| = k-1$. By Claim \ref{claim:finter},
for any $x \in P_k$, $\finter(x) \in S$. It follows (by pigeon hole principle) 
that there exists $x, y \in P_k$
such that $\finter(x) = \finter(y) = I_z \in S$. 
By definition of graph $G$, it contains the four cycle $(x, y, \core_k(x), \core_k(y), x)$.
Since $E(I_z) \supseteq E(G)$, the same four cycle is present in $I_z$ also.
But
by the definition of $\finter(\cdot)$,
$(x, \core_k(x)), (y, \core_k(y)) \notin I_z$. 
Thus it follows that the above four cycle is chordless in $I_z$, which
is a contradiction since $I_z$ is an interval graph.
It follows that $\boxi(G) \ge t - \floor{\sqrt{t}}$.
\end{pf}

\begin{lem} 
\label{lemma:connected-component}
If a multi-graph $M$ has $n$ nodes and at most $n - \frac{n}{r}$ edges, for some $r
\ge 1$, then there is a connected component $C$ in $M$ that has $k$ nodes and
exactly $k-1$ edges, for some $k$ with $1 \le k \le r$. 
\end{lem}

\begin{pf}
Consider those connected components in $M$ where each of them have at least
$r+1$ nodes. 
Call them `large' 
connected components. A connected component which is not large is called 'small'
component. 
Let $g$ and $h$ respectively denote the number of large and small connected components.
Let $n_1, \ldots,
n_g$ respectively be the number of nodes in each of these $g$ large connected
components. Let $n_L = n_1 + \cdots + n_g$ denote the total number of nodes in
the $g$ large components.
Let the total number of edges in these $g$ connected
components together be denoted as $m_L$. Observe that  $m_L \ge n_L - g$.
Let $n_S$ and $m_S$ respectively be the total number of nodes and edges in
the  $h$ small
connected components. 
We have,
$$
n_S = n - n_L ~~~~\mbox{and}~~~~ m_S \le n - \frac{n}{r} - m_L
$$

Recalling that $m_L \ge n_L - g$, we get ~$n_S  \ge m_S + \frac{n}{r} - g$.
Since each large component contains at least $r+1$ nodes, we have
$g < \frac{n}{r}$\,. It follows that $n_S > m_S$.
Thus we can infer that there exists at least one small connected component $C$,
on $k$ nodes and exactly $k-1$ edges, as required.
\end{pf}

\section{Consequences on  Special Classes of Graphs}

\subsection {Chordal Graphs} \label{sec:chordal}

 Let $C$ be a cycle in a  graph $G$. A chord of $C$ is an edge of $G$
 joining two nodes of $C$ which are not consecutive. 
 A graph $G$ is called  a chordal (rigid circuit or triangulated) graph \emph{if and only if}
 every cycle in $G$, of length 4 or more has a chord. Chordal
 graphs arise in many applications (see \cite {Golu}). 
 Chordal graphs constitute one of the most important subclasses of
 perfect graphs \cite {Golu}. 
It is easy to see that the class of interval graphs do not have chordless cycles
of length more than 3, and thus is a subclass of chordal graphs.
Thus it is natural to study the boxicity of chordal graphs.
In fact, the question of bounding the boxicity of a subclass of chordal graphs,
namely the split graphs, was already addressed by Cozzens and Roberts \cite{CozRob}. 
A graph $G$ is a split graph if and only if $G$ and its complement 
$\overline G$ is chordal.  Every 
split graph has the following special structure: Its node set can 
be partitioned into two sets $V_1$ and $V_2$ such that $V_1$ induces
a complete graph in $G$, and $V_2$ induces an independent set. (See 
\cite {Golu}, Chapter 6  for more information on split graphs.) 
Cozzens and Roberts prove the following Theorem.

\begin{thm}[Cozzens and Roberts \cite{CozRob}] \label{thm:CozRob}
If $G$ is a split graph with clique number  $\omega(G)$, then
$\boxi(G) \le \lceil \frac{1}{2}\omega(G) \rceil$.
\end{thm}


The class of split graphs indeed has a very special structure, and is hardly representative
of the much wider class of chordal graphs.
Using our upper bound (Theorem \ref{thm:tw-upper}), we derive an upper bound
for the boxicity of chordal graphs.

The following result is well-known. (See \cite {Diest}, Chapter 12.)

\begin {lem} \label{lemma:chordal}
For a chordal graph $G$,  $\tw(G) = \omega(G) -1$.
\end {lem} 

Combining Lemma \ref {lemma:chordal} with Theorem \ref {thm:tw-upper}, and
noting that $\omega(G) - 1 \le \Delta(G)$ for any graph $G$, where $\Delta(G)$ is its
maximum degree,  we get the following. 

\begin {thm} \label{thm:chordal1}
For a chordal graph $G$, $\boxi(G) \,\le\, \omega(G) + 1 \,\le\, \Delta(G) + 2$.
\end {thm}

\noindent {\bf Sharpness of Theorem \ref {thm:chordal1}}
In \cite{CozRob}, Cozzens and Roberts show that for any $k \ge 1$, there exists a split
graph $G$ such that $\omega(G) = k$ and $\boxi(G) = \lceil \frac{\omega(G)}{2} \rceil$.
Since the class of split graphs is subclass of chordal graphs, it follows that the upper bound of Theorem \ref{thm:chordal1}, is tight up to a factor of $2$.

\subsection {Circular Arc graphs} 
\label{sec:circarc}

A graph $G$ is a circular arc graph if and only if there exists a
one-to-one correspondence between its nodes and a set of arcs of 
a circle, such that two nodes are adjacent \emph{if and only if} the corresponding
arcs intersect. 
Since  the definition of  circular arc  graphs look very
similar to that of interval graphs, it is natural to think that 
their boxicity will be small (possibly bounded above by a
constant), since interval graphs are boxicity 1 graphs. 
In his pioneering paper \cite{Roberts}, F.S. Roberts considers the following graph.

\begin{defn}[Roberts Graph]\label{def:rob}
The Roberts graph on $2n$ nodes is obtained by removing the edges of a 
perfect matching from a complete
graph on $2n$ nodes, 
\end{defn}

Roberts showed that the boxicity of Roberts graph on $2n$ nodes is $n$. A little
inspection will convince the reader that Roberts graph
is indeed a circular arc graph. Thus, 
there exists a circular arc graph of $2n$ nodes, whose boxicity is $n$.

But still, it is possible to get an upper bound for the boxicity of 
circular arc graphs in terms of clique number $\omega(G)$ and thus in terms
of its maximum degree $\Delta(G)$ as follows.

We claim that the pathwidth (and hence treewidth) of a circular
arc graph $G$ is at most $2 \omega(G) -1$. Consider a representation
of $G$ as the intersection graph of  arcs of a circle. 
Let $p_0, \ldots, p_k$, be the end points (left or right) 
of the arcs on this circle, 
as we traverse the circle in the clock-wise direction,
starting from an arbitrarily fixed position.
Let $X_i$ denote the set of nodes of $G$.
whose arcs contain $p_i$. Clearly, $X_i$ for any $i$ induces a complete graph in $G$ and
thus $|X_i| \le \omega(G)$. It is straightforward to verify that
 $(X_1 \cup X_0), \ldots, (X_k \cup X_0)$
constitutes a valid path decomposition of $G$ and thus the pathwidth of $G$ is at most
$2 \omega(G) -1$.

It follows from 
Theorem \ref {thm:tw-upper} that

\begin {thm} \label{thm:circular1}
For a circular arc graph $G$, $$\boxi(G) ~\le~ 2\omega(G) + 1 ~\le~ 2 \Delta(G) + 3.$$
\end {thm} 

\noindent {\bf Tightness of Theorem \ref {thm:circular1}:}
Recall that Roberts graph $G$ on $2n$ nodes (Definition \ref{def:rob}) is a circular arc graph
and its boxicity is $n$.
It is easy to see that $\omega(G) = n$. 
Thus the upper bounds given
by Theorem \ref{thm:circular1}, in terms of $\omega(G)$ and $\Delta(G)$, 
is tight up to a factor $2$ and $4$ respectively.

\subsection {Asteroidal Triple-free Graphs, Co--comparability Graphs and Permutation Graphs}  

An independent set of three nodes in $G$, such that each pair is joined
by a path that avoids the neighborhood of the third is called an
Asteroidal Triple (AT).  A graph is AT-free if and only if it contains
no asteroidal triples. The concept of Asteroidal triples and AT-free graphs
was introduced by Lekkerkerker and Boland \cite {Lekker}, to characterize the 
chordal graphs which are not interval graphs.  They showed that 
a graph  $G$ is an interval graph if and only if, it is simultaneously
a chordal graph and an AT-free graph. 

AT-free graphs generalize (in addition to the class of 
interval graphs) some very important and practically useful
classes of graphs- for example  co-comparability graphs, trapezoidal graphs,
and the permutation graphs. (See \cite {Corn} for a discussion of how 
AT-free graphs are in some sense a unifying generalization of these 
graph classes.) 

In this section we give an upper bound for the boxicity of AT-free graphs
in terms of their maximum degree.  

A \emph{caterpillar} is a tree such that a path (called the \emph{spine})
is obtained by
when all its leaves are deleted. 
In the proof of Theorem 3.16 of \cite{Klok6}, Kloks et al. show that
every connected AT-free graph $G$ has a spanning caterpillar subgraph $T$,
such that adjacent nodes in $G$ are at distance at most four in $T$. Moreover,
for any edge $(u, v) \in E(G)$ with $u$ and $v$ at distance exactly four in $T$,
both $u$ and $v$ are leaves of $T$. Let $p_0, \ldots, p_k$ be the 
nodes along the spline of $G$. Let $X_i$ be the union of $p_i$ and the leaf nodes
attached to $p_i$ in the caterpillar. 
Now it is easy to check that $(X_0\cup X_1 \cup X_2), \ldots, (X_i \cup X_{i+1} \cup X_{i+2}),
\ldots, (X_{k-2} \cup X_{k-1} \cup X_k)$ constitute a path decomposition (and thus
a tree decomposition) of $G$.

\begin {lem}
Let $G$ be an AT-free graph. Then $\tw(G) \le 3\Delta(G) -2$, where $\Delta(G)$
is the maximum degree of $G$.
\end {lem}

\begin{thm}
\label{thm:atfree}
For an AT-free graph $G$,  $\boxi(G) \le 3\Delta(G)$.
\end{thm}

We get better upper bounds when we restrict ourselves to sub classes of AT-free graphs.
Consider the class of 
co--comparability graphs: A graph $G$ is a co- comparability graph if 
and only if its complement is a comparability graph. (See Chapter 5
of \cite {Golu}, for more information on comparability graphs.) 
An interesting characterization of co--comparability graphs is that
they are exactly the class of intersection graphs 
of function diagrams \cite{Gol1}. (A function diagram is a set of curves $\mathbb{C}$,
where each $c_i \in \mathbb{C}$ is a  curve $\{(x, f_i(x)) : 0 \le x \le 1 \}$ for some
$f_i : [0, 1] \rightarrow \mathbb{R}$.)

It is known that 
co-comparability graphs
are properly contained in the class of AT--free graphs, but in turn
is a strict super class of permutation graphs, trapezoidal graphs etc.

\begin {lem}\label{lemma:cococo}
For a co-comparability graph $G$, $\tw(G) \le 2 \Delta(G) - 1$. 
\end {lem} 

\begin{pf}
Let $E(G)$ and $E(\overline{G})$ denote the edge set of $G$ and its complement $\overline{G}$
respectively
and let $V$ be the node set. Let $|V| = n$.
Since $\overline{G}$ is a comparability graph, there exists a partial order $\prec$
in $\overline{G}$ on the node set $V$ such that $(u, v) \in E(\overline{G})$ \emph{if and only if}
$u$ and $v$ are comparable (that is either $u \prec v$ or $v \prec u$).
This partial order gives an orientation to the edge set $E(\overline{G})$, namely,
if $u \prec v$, then the edge $(u, v)$ is directed from $u$ to $v$ and we denote this
directed edge as $[u, v]$.
Define an ordering (i.e., a bijection) $f: V \rightarrow \{1, \ldots, n\}$ for $V$
such that if $(u, v) \in E(\overline{G})$ then $u \prec v$ \emph{if and only if} $f(u) < f(v)$. 
Clearly such an ordering exists for $\overline{G}$;
for instance,
a topological sort on $\overline{G}$ after orienting its
edges as described above, gives such an ordering.
Let $(u, v) \in E(G)$
and $w$ be such that $f(u) < f(w) < f(v)$.  
We claim that
$w$ is adjacent to either $u$ or $v$ or both in $G$.
Assume otherwise. That is, $(u, w), (w, v) \in E(\overline{G})$. 
Since it is given that $f(u) < f(w) < f(v)$,
it follows that
$u \prec w \prec v$ in $\overline{G}$ by the definition of $f(\cdot)$.
Thus by transitivity of $\prec$,  $u \prec v$ and $(u, v) \in E(\overline{G})$,
which is a contradiction.
Having shown that if 
$f(u) < f(w) < f(v)$ then $w$ is adjacent to either $u$ or $v$ or both in $G$,
it is easy to infer that if edge $(u, v) \in E(G)$, then there can be at most $2\Delta(G) - 2$ nodes 
whose $f(\cdot)$ values are between $f(u)$ and $f(v)$. Therefore, 
$|f(u) - f(v)| \le 2\Delta(G) - 1$.
Now 
it is easy to verify that there is a path decomposition for $G$ (and hence a tree decomposition)
$(\{X_i : i \in I\}, T)$, where $I = \{1, \ldots, n\}$ and $T$ is a simple path
$(1, 2, \ldots, n)$, such that 
$X_i = \{u ~|~ i \le f(u) \le i+2\Delta(G) - 1 \}$. It is straightforward to 
verify that the pathwidth of the this path decomposition is $2 \Delta(G) -1$.
\end{pf}

Now it follows from our upper bound (Theorem \ref{thm:tw-upper})
that 

\begin {thm}
\label{thm:coco}
For a co--comparability graph $G$, $\boxi(G) \le 2 \Delta(G) + 1.$
\end {thm}

Permutation graphs are defined as follows.
Let $\pi$ be a permutation of the numbers $1, 2, \ldots, n$, and let
graph $G[\pi] = (V, E)$ be defined as follows: $V = \{1, \ldots, n\}$
and $(i, j) \in E(G)$ \emph{if and only if} $(i - j)(\pi^{-1}(i) - \pi^{-1}(j)) < 0$.
A graph $G$ is a permutation graph if there exists a permutation $\pi$
such that $G$ is isomorphic to $G[\pi]$. It is well-known that the 
complement of a permutation graph is also a permutation graph
 (See Chapter 7, \cite{Golu}).

Permutation graphs are sub-classes of co--comparability graphs 
 (See Chapter 7, \cite{Golu}).
Therefore the above upper bound on boxicity holds also for permutation graphs.
\begin {thm}
\label{thm:perm}
For a permutation graph $G$, $\boxi(G) \le 2 \Delta(G) + 1.$
\end {thm}

\noindent {\bf Tightness of Theorems \ref {thm:atfree}, \ref{thm:coco} and \ref{thm:perm}:}
It is not difficult to see that  Roberts graph  on $2n$ nodes
(Definition \ref{def:rob}) 
is a permutation graph (because its complement  is trivially a permutation graph), 
and hence it is both AT-free and co- comparability.
This graph has maximum degree $2n -1$ and boxicity $n$. 
It follows that  Theorem \ref{thm:coco} and Theorem \ref{thm:perm} are tight up to a factor of $4$
and Theorem \ref{thm:atfree} is tight up to a factor of $6$.

\section {Planar Graph Minors and Boxicity}

The following theorem is well-known.
\begin{thm}[Robertson and Seymour \cite{RS13}] \label{thm:rob-planar}
For every planar graph $H$, there is a constant $c(H)$ such that every graph
with treewidth $\ge c(H)$ has a minor isomorphic to $H$.
\end{thm}

Combining the upper bound theorem (Theorem \ref{thm:tw-upper}) with Theorem \ref{thm:rob-planar},
 we obtain the following. 

\begin{thm}\label{thm:planar-minor}
For every planar graph $H$, there is a constant $c(H)$ such that every graph
with boxicity $\ge c(H)$ has a minor isomorphic to $H$.
\end{thm}

For instance, consider the cycle graph. Note that if a cycle graph on $k$ nodes is a minor of $G$,
then it is also a subgraph of $G$. In other words, by Theorem \ref{thm:planar-minor}, there
exists a constant $c(k)$ such that,
if $\boxi(G) \ge c(k)$ then $G$ contains a cycle on $k$ nodes as a subgraph. 
In the next section,
we show that $c(k)$ is not more than $2k$ by means of a direct approach.

A restatement of Theorem \ref{thm:planar-minor} is as follows.
\begin{thm}
For every planar graph $H$, there is a constant $c(H)$ such that
any minor closed family of graphs which excludes $H$ has boxicity at most $c(H)$.
\end{thm}

\section {Cycles and Boxicity}

The study of various kinds of circuits in graphs is a
well-established area in graph theory.
There is indeed an extensive literature on this topic. 
(Chapter 1 of \cite {Bondy_chap}
gives an introductory survey. See the book by Voss \cite{Voss} for
an extensive treatment.)

\ignore{
	Probably the initial motivation for exploring the properties which causes 
	long simple circuits in graphs was the search for necessary an sufficient 
	condition for hamiltonicity. Many results of  this type are known- the 
	results Egawa et.al. \cite {} and  ... which generalises the classical
	result of Dirac are examples.  
}
A recent result of Birmele \cite{Birmele} relates the 
treewidth of a graph with the length of the longest simple cycle in $G$.
The length of the longest simple cycle in a graph is also known as its
circumference.

\begin {lem}[E. Birmele \cite{Birmele}]\label{lemma:birmele}
For any graph $G$, its circumference is at least $\tw(G) -1$. 
\end {lem}

Combining Lemma \ref {lemma:birmele} with our upper bound Theorem \ref {thm:tw-upper}, we get

\begin {thm} \label{thm:circ}
In any  graph $G$ of boxicity $b$, there exists a simple cycle of length at least $b - 3$.
\end {thm}

\noindent {\bf Sharpness of Theorem \ref{thm:circ}:} It is not difficult to see that
Roberts graph on $2n$ nodes is Hamiltonian. Since its boxicity is known to be $n$,
 the upper bound for Theorem \ref{thm:circ} is tight up to a factor of $2$.

Now we consider induced cycles in a graph.
A cycle in a graph is called
an induced cycle or a chordless cycle if there are no chords for that cycle.
An induced cycle of length $4$ or more is sometimes referred to as a hole,
especially in the perfect graph literature: For example, the strong perfect 
graph theorem states that a graph is perfect if and only if it does not
contain any odd holes.  
The length of the largest induced cycle in a graph $G$ is called the chordality
of $G$.  It may be noted that the chordality of a chordal graph is  $3$, 
that of co-comparability graphs is at most $4$
and 
that of AT--free graphs is at most $5$.
See \cite{CSL} for a survey of graph classes with low chordality.
A great deal of research is done also with respect to 
the existence of cycles in a graph with a given number of chords (diagonals) \cite{Voss}.

The following Lemma follows from Theorem 14  of
\cite {BodThil}.

\begin {lem}[Bodlaender and Thilikos \cite{BodThil}]\label{lemma:bodthil}
Let $G$ be a graph with maximum degree $\Delta$ and chordality $c$.
Then $\tw(G) ~\le~ \Delta^{c -2}$.
\end {lem}

By combining Lemma \ref{lemma:bodthil} with our upper bound  (Theorem \ref {thm:tw-upper}), 
we get the following result

\begin {thm}
Let $G$ be a graph with maximum degree $\Delta$ and boxicity $b \ge 3$.
Then there exists a induced cycle (chordless cycle) of length at least
$\lfloor \log_{\Delta}(b - 2) \rfloor + 2$.
\end {thm} 

\section{Boxicity, Vertex Cover and Related Parameters}

The subset $S \subseteq V(G)$ is called a \emph{vertex cover} of $G$
if every edge of $G$ is incident on at least one vertex from $S$.
A vertex cover of minimum cardinality is called a \emph{minimum vertex
cover}. We denote the cardinality of a minimum vertex cover of $G$  by 
$M\!V\!C(G)$. It is easy to observe that if $S$ is a minimum vertex 
cover of $G$, then $V(G) - S$ forms an independent set of $G$. Thus,
$M\!V\!C(G) = n - \alpha(G)$, where $\alpha(G)$ is the independence number
(the cardinality of the maximum independent set) of $G$. It is easy to prove the following
Lemma. 

\begin {lem}
For any graph $G$, $\tw(G) \le M\!V\!C(G)$. 
\end {lem}

Now, applying Theorem \ref{thm:tw-upper} we get:

\begin {thm}
\label {thm:mvc-boxicity}
For any graph $G$, $\boxi(G) \le M\!V\!C(G) + 2 = n - \alpha(G) + 2$.
\end {thm}

It is interesting to investigate whether the above bound in terms of $M\!V\!C(G)$ can
be further tightened. For instance, we can show that
if a graph $G$ has a vertex cover which induces a complete graph, then
$\boxi(G) \le \lceil (M\!V\!C(G)+1)/2 \rceil$. 
To see this, first we recall from Section \ref{sec:chordal} that
if the node set $V$ of a graph $G$ can be partitioned into $V_1$ and $V_2$
such that $V_1$ induces a complete graph and $V_2$ induces an independent
set, then $G$ is a split graph. It follows that
if a vertex cover in $G$ induces a complete graph, then $G$ is a split graph.
Now the above bound follows from Theorem \ref{thm:CozRob}, because
$\omega(G) -1 \le M\!V\!C(G) \le \omega(G)$ for a split graph $G$.


A set of \emph{dominating edges} $D$ is a collection of edges of $G$ such
that  any edge in $E(G) - D$ is  adjacent to at least one edge in $D$. 
For example the reader may notice that any maximal matching in $G$
constitutes a dominating edge set.
A dominating edge set of minimum cardinality is called a minimum 
dominating edge set. We denote the cardinality of the minimum dominating edge set
of $G$ by $M\!E\!D(G)$.
It is easy to see that 
$M\!V\!C(G) \le 2\, M\!E\!D(G)$.
Combining this with Theorem \ref {thm:mvc-boxicity}, we have

\begin {thm}
\label {thm:dominating-boxicity}
For any graph $G$,  $boxi(G) \le 2\, M\!E\!D(G) + 2$. 
\end {thm} 

In this connection, we note that Cozzens and Roberts \cite{CozRob}
had proved the following: for any graph $G$, $\boxi(G) \le M\!E\!D(\overline{G})$.
Clearly our result (Theorem
\ref {thm:dominating-boxicity}) complements their result, by showing that
$M\!E\!D(G)$ itself can control the boxicity
of $G$. Thus,

\begin {thm}
\label {thm:dominating-boxicity2}
For any graph $G$, $\boxi(G) \le \min\{2\, M\!E\!D(G) + 2, M\!E\!D(\overline{G})\}.$
\end {thm} 

Now, let us consider a variant of minimum vertex cover, namely the 
\emph{minimum feedback vertex cover}. A feedback vertex cover
$S$ is a subset of $V(G)$ such that the induced subgraph on $V-S$ 
is a forest. A feedback vertex cover of minimum cardinality is
called a minimum feedback vertex cover, and we denote its cardinality
by $M\!F\!V\!C(G)$. Clearly every vertex cover of $G$ is a feedback vertex
cover also, and thus $M\!F\!V\!C(G) \le M\!V\!C(G)$.  The reader may also
note that in general $M\!F\!V\!C(G)$ can be much smaller than $M\!V\!C(G)$: For
example for a cycle on $n$ nodes, $\frac {M\!V\!C(G)}{M\!F\!V\!C(G)} = \Omega(n)$.

\begin {thm} \label{thm:mfvc}
$\boxi(G) \le M\!F\!V\!C(G) + 3.$
\end {thm}

\begin {pf}
Let $S$ be a {minimum feedback vertex cover} of $G$. Since the induced
subgraph on $V-S$ is a forest,  there exists a tree decomposition
$(\{X_i : i \in I\}, T),$ whose width is $1$. Clearly, $(\{X_i \cup S : i \in I\}, T)$
is a valid tree  decomposition of $G$ whose width is $|S| + 1 = M\!F\!V\!C(G) + 1$.
Now, applying the upper bound theorem (Theorem \ref{thm:tw-upper}), the result follows.
\end {pf}

\noindent
\textbf{Sharpness of Theorems \ref{thm:dominating-boxicity} and \ref{thm:mfvc}:}
It is easy to check that for Roberts' graph  $G$ on $2n$ nodes (see Definition \ref{def:rob}), 
$M\!V\!C(G) = 2n -2$,
whereas its boxicity  is $n$. Thus, the upper bound of Theorem \ref{thm:dominating-boxicity} is tight up to a factor of $2$. Similarly, it is easy to see that $M\!F\!V\!C(G) \ge 2n - 4$, and thus 
the upper bound of Theorem \ref{thm:mfvc} is also tight up to a factor of $2$.

\section{Algorithmic Consequence}

\begin{thm}
For a  bounded treewidth graph $G=(V,E)$ on $n$ nodes,
a box representation of $G$ in constant dimension can be constructed
in $O(n)$ time.
\end{thm}

\begin{pf}
We construct the interval graph representation of $G$ using $\tw(G) + 2$ interval graphs
$I_0, \ldots, I_{\tw(G) + 1}$, as described in the proof of Theorem \ref{thm:tw-upper}.
It is not difficult to observe that the proof of Theorem \ref{thm:tw-upper} is constructive.
It remains to show that this construction can be implemented in linear time when
$\tw(G)$ is bounded.
(Recall from Section \ref{sec:boxi-inter} that 
the interval graph representation of $G$ is equivalent to its box representation.)

It is well-known (see for instance \cite{Bodland2}) that
if $\tw(G) \le k$ then $|E(G)| \le k n - \frac{1}{2} k (k+1)$.

We convert the constructive proof of Theorem \ref{thm:tw-upper} into a linear time
algorithm consisting of the following steps.
We show that each of these steps can be implemented in linear time.

\begin{enumerate}
\item 
Given a  bounded treewidth graph $G$, Bodlaender \cite{Bodland2} gives
an $O(n)$ algorithm to construct the optimum tree decomposition 
 $(\{X_i : i \in I\}, T)$ of $G$. (In this tree decomposition, $|I| = O(n)$.)

\item \label{step2}
Convert this tree decomposition into a normalized tree  decomposition
$(\{X_i : i \in I\}, r,  T)$ as described in the proof of Lemma \ref{lemma:norm-tree}.
It is straightforward to verify that this conversion takes $O(n)$ time.
It is also easy to see that, while doing this conversion, we can additionally 
obtain the following,
without increasing the time complexity:
\begin{enumerate}
\item
 An ordering of $I$,
sorted in the non-increasing order of their heights in the rooted tree $T$.
\item
The $b(u)$ and $\level(u)$ values for each $u \in V$.
\end{enumerate}

\item
Compute $\theta(u)$ for each $u \in V$ as described in the proof of Lemma \ref{lemma:theta}.
The sorted order of $I$ as required in this proof is already computed in Step \ref{step2}.
The remaining steps in this proof can be implemented in a straightforward way such
that it takes only constant time for computing $\theta(u)$ for each node $u \in V$.
Thus the total time is $O(n)$.

\item
Now we construct the interval graphs $I_0, \ldots, I_{\tw(G)}$. 
To construct $I_i$, $0 \le i \le \tw(G)$, we need to compute $\ell_i(v)$ and $r_i(v)$ for
each $v \in V$, as described in Section \ref{sec:up}.
If $\theta(v) = i$ then computing $\ell_i(v)$ and $r_i(v)$ is trivial.
If $\theta(v) \not= i$, then first we have to compute the set $S = \theta^{-1}(i) \cap N(v)$.
It is easy to see that this can be done in $O(|N(v)|)$ time. After this, we have
to compute $\min_{u \in S} r_i(u)$. This takes additional $O(|S|)$ time.
Thus the total time taken for node $v$ is $O(|N(v)| + |S|) = O(|N(v)|)$.
Hence the overall time taken to construct $I_i$ is $O(|E(G)|) = O(n)$. (Recall that
for a bounded treewidth graph $G$, $O(|E(G)|) = O(n)$.)

\item
It remains to construct the interval graph $I_{\tw(G) + 1}$ as described
in Section \ref{sec:up}. For this, it is required to compute the depth first traversal order
of $T$. It is trivial to see that such an ordering can be computed in $O(|I|) = O(n)$ time.
Assigning for each $u \in V$, its corresponding sequence numbers in the traversal order
  to $\ell_i(u)$ and $r_i(u)$
can be done without increasing the time complexity.
\end{enumerate}
\end{pf}

\end {document}